\documentclass[11pt,a4paper]{article}
\usepackage[T1]{fontenc}
\usepackage{lmodern,microtype}
\usepackage[margin=1in]{geometry}
\usepackage{amsmath,amssymb,amsthm,mathtools}
\usepackage[hidelinks,pdfusetitle]{hyperref}
\allowdisplaybreaks[1]
\newtheorem{theorem}{Theorem}
\newtheorem{lemma}[theorem]{Lemma}
\newtheorem{corollary}[theorem]{Corollary}
\newcommand{\E}{\mathbb E}
\newcommand{\PP}{\mathbb P}
\newcommand{\R}{\mathbb R}
\newcommand{\Var}{\operatorname{Var}}
\newcommand{\sd}{\operatorname{sd}}
\newcommand{\dK}{d_{\mathrm K}}
\newcommand{\ind}{\mathbf 1}
\newcommand{\pairs}{\mathcal P}

\title{A truncated-moment Berry--Esseen bound for dependency graphs}
\author{Maximilian Janisch\thanks{\href{mailto:mail@maximilianjanisch.com}{\nolinkurl{mail@maximilianjanisch.com}}}}
\date{26 November 2026}

\begin{document}
\maketitle

\begin{abstract}
Let $X_1,\ldots,X_n$ be centered, square-integrable random variables with
a dependency graph of maximum degree $D$. Put $S=\sum_i X_i$ and suppose
that $\sigma^2=\Var S>0$. We prove that the Kolmogorov distance between
$S/\sigma$ and the standard normal distribution is at most
\[
 C\,\frac{D+1}{\sigma^2}\sum_{i=1}^n
 \E\!\left[X_i^2\min\!\left\{1,\frac{(D+1)|X_i|}{\sigma}\right\}\right].
\]
For $2<\delta\le3$, this gives the bound
$C(D+1)^{\delta-1}\sigma^{-\delta}\sum_i\E|X_i|^\delta$ with one
absolute constant $C$ for the entire range. This proves the conjecture of
Janisch and Lehéricy (2024), including the third-moment case conjectured
by Raič in 2003. The proof combines Stein's method, a deterministic
inequality for adjacent pairs, and induction on the number of summands.
A formalization in Lean~4 accompanies the paper.
\end{abstract}

\section{Introduction}
During my PhD, Thomas Lehéricy and I studied normal approximation for sums
of random variables with a dependency graph.
In our paper~\cite[Lemma~5.1]{JL}, we sharpened a cumulant
estimate that my advisor, Ashkan Nikeghbali, had developed with Valentin
Féray and Pierre-Loïc Méliot in their work on mod-$\phi$
convergence~\cite[Theorem~9.8]{FMN}. This gave us bounds on the Kolmogorov
distance between the standardized sum and the standard normal distribution.
For some ranges of the maximum degree $D$, our bounds improved on the
Berry--Esseen bound of Chen and Shao~\cite[Theorem~2.7]{ChenShao}.
For bounded summands our Fourier argument gave an essentially optimal
estimate. With only a moment of order $2<\delta\le 3$, it gave a slower
rate~\cite[Theorems~2.1 and~2.4 and Section~4]{JL}.

A simple example motivated our conjecture. Suppose the graph consists of
disjoint cliques of size $D+1$, with identical summands in each clique.
Distinct cliques are independent by the definition of a dependency graph.
After grouping, $N$ summands give just $N/(D+1)$ independent terms. We
conjectured that the Berry--Esseen estimate for independent random variables
would still hold for every dependency graph of maximum degree $D$, with
$N/(D+1)$ in place of $N$ in the rate~\cite[Conjecture~1.4]{JL}.

Martin Raič had already conjectured the case $\delta=3$ in
2003~\cite[Conjecture~5.4]{Raic}. Just before stating his conjecture, he
suggested that the induction used in the independent case might also work
for local dependence. For $m$-dependent sequences, Shergin had proved the
corresponding moment bound in 1979~\cite[Theorem~2]{Shergin}. For general
dependency graphs, Cai, Shao and Zhang~\cite[Theorem~3.1]{CSZ} recently
obtained a bound under a fourth-moment assumption.

I finished my PhD without a proof of the conjecture. As AI models improved,
I kept giving them the problem. I first obtained a complete proof with
OpenAI's GPT-6.2~Astra in November 2026. The complete Lean~4 formalization
and instructions for checking the proof are included in the arXiv ancillary files.

The proof gives a stronger result. Let $W=\sum_i X_i$ have variance one,
where the $X_i$ are centered and square-integrable and have a dependency
graph of maximum degree $D$. Write $q=D+1$, let $Z$ be standard normal,
and denote Kolmogorov distance by $\dK$. Then
\[
 \dK(W,Z)\le 250{,}000\,q\sum_i
 \E\!\left[X_i^2\min\{1,q|X_i|\}\right].
\]
At $q=1$ this recovers the truncated-moment form of Osipov's inequality
for independent summands~\cite{Osipov}. Since $t^2\min\{1,t\}\le t^\delta$
for every $t\ge0$ and $2<\delta\le3$, the same constant gives
$\dK(W,Z)\le 250{,}000\,q^{\delta-1}\sum_i\E|X_i|^\delta$ throughout
the conjectured range. Theorem~\ref{thm:published} states the result in
the notation of our original conjecture.

\paragraph{Idea of the proof.}
We use Stein's method~\cite{Stein72} and induct on the number of summands,
as in the independent case~\cite{Bolthausen,Raic}. Smoothing the indicator
of a half-line in Stein's equation gives a sum over adjacent pairs
$(i,j)$, including $i=j$. Each term asks whether $W$ lies in a short
interval. The difficulty is that its width depends on all the summands
in the neighborhoods of $i$ and $j$. Chen and Shao estimated interval
probabilities under local dependence by concentration inequalities~\cite[Section~3]{ChenShao}.

Lemma~\ref{lem:charge} is the step that makes induction possible.
Its proof uses the fact that neighbors carrying at least half the average
weight in a closed neighborhood account for at least half its total weight.
After summing over pairs, the lemma gives intervals whose
widths and weights depend only on $X_i$ and $X_j$. The sum tested against
each interval is now $W_{ij}$, obtained by deleting both closed
neighborhoods. The dependency graph makes $W_{ij}$ independent of
$(X_i,X_j)$. We can therefore fix the two endpoint values and use the
induction hypothesis on the smaller graph to estimate the interval
probability. Each vertex belongs to at most $q$ closed neighborhoods.
Summing the resulting estimates gives the truncated moment above.

\section{Main results}\label{sec:results}
Let $G$ be a simple graph on a finite set $V$. We call $G$ a
\emph{dependency graph} for $(X_i)_{i\in V}$ if the random vectors
$(X_i)_{i\in A}$ and $(X_j)_{j\in B}$ are independent whenever $A,B$
are disjoint vertex sets with no edge between them. Write $A_i$ for the
closed neighborhood of $i$, consisting of $i$ and its neighbors.
Fix a positive integer $q\ge\max_i|A_i|$. In particular, $q=D+1$ is
admissible if $D$ is the maximum degree. Let $Z$ be standard normal
with distribution function $\Phi$. For real random variables $U,V'$,
write $\dK(U,V')=\sup_z|\PP(U\le z)-\PP(V'\le z)|$.

\begin{theorem}[Truncated-moment bound]\label{thm:main}
Let $(X_i)_{i\in V}$ be centered and square-integrable with dependency
graph $G$. Suppose that $W=\sum_iX_i$ has variance one. Put
\begin{equation}\label{eq:ell}
 \ell(t)=t^2\min(1,t)\quad(t\ge0),\qquad
 \beta=\frac1q\sum_i\E\ell(q|X_i|).
\end{equation}
Then $\dK(W,Z)\le\min(1,250{,}000\,\beta)$.
\end{theorem}

\begin{corollary}[Centered moment bound]\label{cor:jl}
Under the hypotheses of Theorem~\ref{thm:main}, if $2<\delta\le3$ and
each $X_i$ has a finite absolute moment of order $\delta$, then
\[
 \dK(W,Z)\le\min\!\left\{1,\;250{,}000\,q^{\delta-1}\sum_i\E|X_i|^\delta\right\}.
\]
\end{corollary}
\begin{proof}
For $0\le t\le1$ we have $t^3\le t^\delta$, and for $t\ge1$ we have
$t^2\le t^\delta$. Hence $\ell(t)\le t^\delta$ and
$\beta\le q^{-1}\sum_i\E(q|X_i|)^\delta$.
\end{proof}

Taking $q=D+1$ and $\delta=3$ gives Raič's Conjecture~5.4~\cite{Raic}.
Centering and rescaling now give Conjecture~1.4 of~\cite{JL} in its
published form.

\begin{theorem}[Published form of the conjecture]\label{thm:published}
Let $(Y_i)_{i\in V}$ be real random variables with dependency graph $G$ of
maximum degree $D$. Fix $2<\delta\le3$ and suppose that every $Y_i$ has finite
$\delta$-th absolute moment. Put $S=\sum_iY_i$, assume $\sigma^2=\Var S>0$, and set
\[
 W=\frac{S-\E S}{\sigma},\qquad N=|V|,\qquad q=D+1,\qquad
 \mathcal A_\delta=\sum_i\E|Y_i|^\delta,\qquad
 \xi_\delta=\left(\frac{N}{\mathcal A_\delta}\right)^{1/\delta}
              \sqrt{\frac{\sigma^2}{Nq}}.
\]
Then for every $z\in\R$,
\[
 \big|\PP(W\le z)-\Phi(z)\big|
 \le \min\!\left\{1,\;2{,}000{,}000\,\xi_\delta^{-\delta}
                   \left(\frac q N\right)^{(\delta-2)/2}\right\}.
\]
\end{theorem}
The quantity $\sigma$ is denoted by $v$ in~\cite{JL}.
The file
\texttt{JanischLehericy/Published.lean} gives the same hypotheses, definitions,
and conclusion in the theorem \texttt{conjecture\_one\_four}. Conjecture~1.3 of
the arXiv version of \cite{JL} allows $\sum_i\E|Y_i-c_i|^\delta$ for arbitrary
constants $c_i$. Applying Theorem~\ref{thm:published} to $(Y_i-c_i)_{i\in V}$
gives this formulation, since translation preserves block independence and the
standardized sum.

\begin{proof}[Proof of Theorem~\ref{thm:published} from Corollary~\ref{cor:jl}]
Set $X_i=(Y_i-\E Y_i)/\sigma$. Their sum is $W$. By the triangle inequality
in $L^\delta$, $\|Y_i-\E Y_i\|_\delta\le2\|Y_i\|_\delta$. Since $2^\delta\le8$,
replacing centered moments by raw moments costs at most a factor of eight.
Substituting the definition of $\xi_\delta$ gives
\[
 \xi_\delta^{-\delta}(q/N)^{(\delta-2)/2}
 = q^{\delta-1}\sigma^{-\delta}\mathcal A_\delta .
\]
The constant in Theorem~\ref{thm:published} is therefore $8\cdot250{,}000$.
\end{proof}

\section{An inequality for adjacent pairs}\label{sec:pairs}
Fix real numbers $v_i$ on the vertices of $G$ and bounds $a_i\ge |v_i|$.
We sum over ordered pairs $\pairs=\{(i,j):j\in A_i\}$, including $(i,i)$.
Write
\begin{equation}\label{eq:pair-notation}
 w=\sum_i v_i,\quad R_i=\sum_{k\in A_i}a_k,\quad
 A_{ij}=A_i\cup A_j,\quad w_{ij}=\sum_{k\notin A_{ij}}v_k,\quad
 L_{ij}=q\max(a_i,a_j).
\end{equation}
Fix a location $u\in\R$ and a smoothing width $\alpha>0$. The function
\begin{equation}\label{eq:F}
 F(s)=\min(1,s/\alpha)\,\ind\{|w-u|\le\alpha+s\},\qquad s\ge0,
\end{equation}
is nondecreasing. Define
\begin{equation}\label{eq:Pu}
 P_u=\sum_{(i,j)\in\pairs}a_i a_j F(R_i+R_j)
 =\sum_{(i,j)\in\pairs}a_i a_j\min\Big(1,\frac{R_i+R_j}{\alpha}\Big)
  \ind\{|w-u|\le\alpha+R_i+R_j\}.
\end{equation}

\begin{lemma}\label{lem:charge}
For every such choice of numbers,
\begin{align}
 P_u
 \le 8\sum_{(i,j)\in\pairs}a_i a_j\min(1,4L_{ij}/\alpha)
       \ind\{|w_{ij}-u|\le\alpha+10L_{ij}\}.       \label{eq:charge}
\end{align}
\end{lemma}
\begin{proof}
Set $J=\sum_i a_iR_i F(2R_i)$. By symmetry of $\pairs$ and monotonicity of $F$, $P_u$ is at most $2J$. For each $i$ retain those $j\in A_i$ for which $q a_j\ge R_i/2$. The discarded $a_j$ sum to at most $R_i/2$, because $|A_i|\le q$. Thus
\[
 J\le2\sum_{(i,j)\ {\rm retained}} a_i a_j F(2R_i).
\]
Split the retained pairs according to the size of $R_j$.

\emph{If $R_j>4L_{ij}$,} retention implies $R_i\le2q a_j$, and $R_j>4L_{ij}$ implies $\max(a_i,a_j)<R_j/(4q)$. Hence $R_i<R_j/2$ and $a_i<R_j/(4q)$. There are at most $q$ incoming $i$ at a fixed $j$, so their $a_i$ sum to at most $R_j/4$. Their total contribution to $J$, including the factor $2$, is bounded by
\[
 2\sum_j a_j F(2R_j)\frac{R_j}{4}=J/2.
\]
\emph{If $R_j\le4L_{ij}$,} retention gives $R_i\le2L_{ij}$. Whenever $F(2R_i)>0$,
\[
 |w_{ij}-u|\le |w-u|+R_i+R_j
 \le \alpha+3R_i+R_j\le\alpha+10L_{ij}.
\]
Also $\min(1,2R_i/\alpha)\le\min(1,4L_{ij}/\alpha)$. We may now sum this upper bound over every ordered pair.
Writing $\Sigma$ for the sum on the right of~\eqref{eq:charge}
without its factor $8$, the two cases give $J\le J/2+2\Sigma$, hence
$J\le4\Sigma$. The lemma follows from $P_u\le2J$.
\end{proof}

\section{From Stein's identity to the pair inequality}\label{sec:stein}
Apply Section~\ref{sec:pairs} with $v_i=X_i$, $w=W$ and
\begin{equation}\label{eq:envelopes}
 a_i=\max(|X_i|,\|X_i\|_2).
\end{equation}
In this setting $P_u$ in~\eqref{eq:Pu} is a random variable and the
deleted sum is $W_{ij}=\sum_{k\notin A_{ij}}X_k$. No vertex outside
$A_{ij}=A_i\cup A_j$ is adjacent to $i$ or $j$. Hence $W_{ij}$ is
independent of $(X_i,X_j)$, and thus of $(a_i,a_j)$.

\begin{lemma}[Stein solution]\label{lem:stein-sol}
Let $z\in\R$, and let $\varphi$ be the standard normal density.
The bounded continuous solution of the normal Stein equation
\begin{equation}\label{eq:stein-eq}
 f_z'(x)-xf_z(x)=\ind\{x\le z\}-\Phi(z),\qquad x\ne z,
\end{equation}
is
\[
 f_z(x)=\frac{\Phi(\min(x,z))(1-\Phi(\max(x,z)))}{\varphi(x)}.
\]
It satisfies $0<f_z\le\sqrt{\pi/2}$. The distributional derivative
$\mu_z$ of $f_z'$ is a signed measure of total variation
$|\mu_z|(\R)\le2$.
\end{lemma}
\begin{proof}
The displayed solution follows by multiplying~\eqref{eq:stein-eq} by
$\varphi(x)$ and integrating from $-\infty$ to $x$. It is positive and at most $(1-\Phi(|x|))/\varphi(|x|)\le\sqrt{\pi/2}$.
Away from $z$, equation~\eqref{eq:stein-eq} gives $f_z'(x)=g_z(x)+\ind\{x\le z\}-\Phi(z)$
with $g_z(x)=xf_z(x)$, so the variation of $f_z'$ is at most two once we show that $g_z$ is increasing from $-(1-\Phi(z))$ to $\Phi(z)$. Indeed, with $r(t)=(1-\Phi(t))/\varphi(t)$, the signs of the derivatives of $g_z$ on $x<z$ and $x>z$ are those of $(1+x^2)r(-x)+x$ and $(1+x^2)r(x)-x$, respectively. They are nonnegative by $r(t)\ge t/(1+t^2)$ for $t\ge0$. For completeness, integration by parts gives
$1-\Phi(t)=\varphi(t)/t-\int_t^\infty\varphi(s)s^{-2}\,ds$ for $t>0$, yielding $t/(1+t^2)\le r(t)\le1/t$ and the claimed limits. Also $r(t)=\int_0^\infty e^{-tv-v^2/2}\,dv\le r(0)$ for $t\ge0$.
\end{proof}

\begin{lemma}\label{lem:stein}
For $0<\alpha\le1$,
\[
 \dK(W,Z)\le 2\alpha+6\sup_{u\in\R}\E P_u .
\]
\end{lemma}
\begin{proof}
\emph{Smoothing.}
Fix $z\in\R$, take $U$ uniform on $[0,\alpha]$ and independent of all
the summands, and put $f(x)=\E_U f_z(x+U)$. Averaging translates of the
continuous, piecewise differentiable function $f_z$ gives
$f'(x)=\E_U f_z'(x+U)$ and makes $f'$ continuous. For each $x$, almost
surely $x+U\ne z$. Averaging~\eqref{eq:stein-eq} at $x+U$ therefore gives
\[
 f'(x)-xf(x)=\PP_U(x+U\le z)-\Phi(z)+\E_U[Uf_z(x+U)] .
\]
This identity holds for every $x$, so we can substitute $x=W$ even at
atoms. The last term costs at most $\alpha\sqrt{\pi/2}/2$ by
Lemma~\ref{lem:stein-sol}. Since
$\PP(W+U\le z)\le\PP(W\le z)\le\PP(W+U\le z+\alpha)$ and
$\Phi(z+\alpha)-\Phi(z)\le\alpha/\sqrt{2\pi}$, we obtain
\[
 \dK(W,Z)\le 2\alpha+\sup_z\big|\E[f'(W)-Wf(W)]\big|,
\]
where $f$ depends on $z$.

\emph{Expansion over pairs.}
Write $X_{A_i}=\sum_{k\in A_i}X_k$ and $c_{ij}=\E[X_iX_j]$. Let $s_i$ be the
secant slope $(f(W)-f(W-X_{A_i}))/X_{A_i}$ when $X_{A_i}\ne0$, and set $s_i=f'(W)$
when $X_{A_i}=0$. Centering and independence of $X_i$ from $W-X_{A_i}$ give
\[
 \E[Wf(W)]=\sum_i\E\big[X_i(f(W)-f(W-X_{A_i}))\big]
           =\sum_{(i,j)\in\pairs}\E[X_iX_j s_i].
\]
Also $\sum_{\pairs}c_{ij}=1$ by $\Var W=1$, since $c_{ij}=0$ for
$j\notin A_i$ by independence and centering. Thus
$\E f'(W)=\sum_{\pairs}c_{ij}\E f'(W)$. Independence of $(X_i,X_j)$ and $W_{ij}$ gives
$\E[(c_{ij}-X_iX_j)f'(W_{ij})]=0$, and therefore the exact identity
\begin{align}
 \E[f'(W)-Wf(W)]={}&
 \sum_{\pairs}\E[(c_{ij}-X_iX_j)(f'(W)-f'(W_{ij}))]\notag\\
 &+\sum_{\pairs}\E[X_iX_j(f'(W)-s_i)] .\label{eq:pair-identity}
\end{align}
For each outcome, the mean value theorem gives $s_i=f'(t)$ at a point $t$
between $W-X_{A_i}$ and $W$. Both $|W-t|$ and $|W-W_{ij}|$ are at most $R_i+R_j$.
Each of $|c_{ij}|$ and $|X_iX_j|$ is at most $a_i a_j$.
Thus the increments of $f'$ in~\eqref{eq:pair-identity} occur over
distances of at most $R_i+R_j$, with total absolute weight at most
$3a_i a_j$ per pair.

\emph{Bounding the increments.}
Choose the right-continuous version of $f_z'$. By Lemma~\ref{lem:stein-sol},
$f_z'(y)=c+\mu_z((-\infty,y])$ for a constant $c$. Averaging gives
$f'(x)=c+\int_\R H(x-u)\,\mu_z(du)$, where $H(y)=\PP(U\ge-y)$
increases linearly from zero to one on $[-\alpha,0]$ and is constant
on either side. Hence $|H(y)-H(y-v)|\le\min(1,|v|/\alpha)$, and this
difference vanishes unless $|y|\le\alpha+|v|$. It follows that
\begin{equation}
 |f'(x)-f'(x-v)|\le\int_{\R}
 \min(1,|v|/\alpha)\ind\{|x-u|\le\alpha+|v|\}\,|\mu_z|(du).
 \label{eq:bv}
\end{equation}
The right side of~\eqref{eq:bv} is nondecreasing in $|v|$. Applying
it with $x=W$ and $|v|\le R_i+R_j$ to every term
of~\eqref{eq:pair-identity}, and recalling the definition of $P_u$,
bounds the absolute value of~\eqref{eq:pair-identity} by $3\int \E P_u\,|\mu_z|(du)\le6\sup_u\E P_u$. This bound is uniform in $z$.
\end{proof}

\section{Moment estimates and induction}\label{sec:induction}
We use the function $\ell$ and parameter $\beta$ from~\eqref{eq:ell}
to bound the cost of deleting the two neighborhoods.

\begin{lemma}\label{lem:mom}
The envelopes and deleted sums satisfy
\begin{align}
 \sum_{(i,j)\in\pairs}\E[a_i a_j\min(1,L_{ij})]&\le 8\beta, \label{eq:moment}\\
 \|W-W_{ij}\|_2^2&\le \tfrac1{16}+16\beta. \label{eq:variance}
\end{align}
\end{lemma}
\begin{proof}
Let $T$ be a nonnegative square-integrable random variable and $s=\|T\|_2$. Then $\E[T^2\ind\{T\ge s/\sqrt{2}\}]\ge s^2/2$. Hence $\E\ell(T)\ge \ell(s)/(2\sqrt2)$, and $\E\ell(\max(T,s))\le(1+2\sqrt2)\E\ell(T)$. Moreover, pointwise,
\[
 a_i a_j\min(1,L_{ij})\le q^{-2}\big(\ell(q a_i)+\ell(q a_j)\big).
\]
Each vertex occurs at most $q$ times in each position. Sum and apply the preceding expectation bound to $T=q|X_i|$ to obtain $2(1+2\sqrt2)\beta\le8\beta$.

For $t\ge0$ we have $t^2\le 1/64+8\ell(t)$. The constant covers $t\le1/8$,
and $\min(1,t)\ge1/8$ covers $t\ge1/8$.
Since $W-W_{ij}=\sum_{k\in A_{ij}}X_k$ and $|A_{ij}|\le2q$,
Cauchy--Schwarz and this scalar inequality give
\[
 \|W-W_{ij}\|_2^2
 \le \frac{2}{q}\sum_{k\in A_{ij}}\E\big[(q|X_k|)^2\big]
 \le \frac{|A_{ij}|}{32q}+16\beta
 \le \frac1{16}+16\beta . \qedhere
\]
\end{proof}

\begin{proof}[Proof of Theorem \ref{thm:main}]
Set $C_0=250{,}000$ and $M=50{,}000$. We prove $\dK(W,Z)\le C_0\beta$ by strong induction on $|V|$, allowing every admissible $q$ at each size.

If $\beta\ge1/M$, the estimate follows from $\dK\le1$. Otherwise take $\alpha=M\beta\in(0,1)$. Positivity of $\alpha$ follows from $\Var W=1$.

\emph{Applying induction to the deleted sums.}
For each pair let $\sigma_{ij}=\sd(W_{ij})$. By \eqref{eq:variance} and the triangle inequality, $\sigma_{ij}\ge1-\sqrt{1/16+16\beta}>1/2$. Thus in this branch the outside sum has positive variance and has strictly fewer vertices
(since $i\in A_{ij}$). The induction starts vacuously, as the empty sum has zero variance. The graph induced by $G$ on $V\setminus A_{ij}$ is a dependency graph for $(X_k)_{k\notin A_{ij}}$, and its closed neighborhoods still have size at most $q$. For the normalized family $(X_k/\sigma_{ij})_{k\notin A_{ij}}$, whose sum is $W_{ij}/\sigma_{ij}$, the quantity defined as $\beta$ in~\eqref{eq:ell} is at most $8\beta$, since $\ell(t/\sigma_{ij})\le\max(1,\sigma_{ij}^{-3})\ell(t)$. Induction gives
\[
 \dK(W_{ij},\sigma_{ij}Z)\le 8C_0\beta .
\]
Consequently for every deterministic $h\ge0$ the closed-interval probability $\PP(|W_{ij}-u|\le h)$ is at most $\min(1,2h)+16C_0\beta$. Here the normal density is at most $1/(\sigma_{ij}\sqrt{2\pi})<1$, and the cdf error enters at both endpoints, with a left limit at the lower endpoint.

\emph{Conditioning on the pair.}
The right side of \eqref{eq:charge} has no restrictions involving $R_i$ or $R_j$. Its weights $a_i,a_j,L_{ij}$ depend only on $(X_i,X_j)$, which is independent of $W_{ij}$. Hence the preceding interval estimate applies conditionally with $h=\alpha+10L_{ij}$. For $L\ge0$ and $b=\min(1,4L/\alpha)$ we have
\[
 b\min(1,2\alpha+20L)\le\min(1,28L)\le28\min(1,L),
 \qquad b\le (4/\alpha)\min(1,L).
\]
The first bound uses $2\alpha b\le8L$ and $20Lb\le20L$.
Together with Lemma~\ref{lem:charge} and~\eqref{eq:moment}, this gives, for every $u\in\R$,
\[
 \E P_u\le 8\,(28+64C_0\beta/\alpha)
 \sum_{(i,j)\in\pairs}\E[a_i a_j\min(1,L_{ij})]
 \le 64\,(28+64C_0\beta/\alpha)\,\beta .
\]
\emph{Closing the induction.} By Lemma~\ref{lem:stein},
\begin{align*}
 \dK(W,Z)
 \le 2\alpha+6\sup_{u\in\R}\E P_u
 \le 2\alpha+10{,}752\,\beta+24{,}576\, C_0\beta^2/\alpha .
\end{align*}
Substituting $\alpha=M\beta$ gives coefficient $2M+10{,}752+24{,}576\,C_0/M=233{,}632<C_0$. This closes the induction.
\end{proof}

\section*{Acknowledgments}
I thank Thomas Lehéricy for fruitful discussions and for his suggestions
on the presentation of the proof.

\section*{Use of AI}
I discussed the conjecture with earlier models, including OpenAI's
GPT-6~Astra and Anthropic's Claude Fable~5.2. GPT-6.2~Astra generated the
Lean~4 proof and drafted the mathematical exposition following the
introduction. Thomas Lehéricy used Claude Opus~5.5 in preparing his
suggested revisions to the exposition.

\end{document}